\documentclass[12pt]{article}
\usepackage{amsfonts,amssymb}
\usepackage{latexsym}
\usepackage{theorem}
\usepackage{amsmath}

\IfFileExists{mathrsfs.sty}{\usepackage{mathrsfs}\let\mathcal\mathscr}{}
\newcommand{\linespacing}{1.25}
\renewcommand{\baselinestretch}{\linespacing}
\newtheorem{theorem}{Theorem}[section]
\newtheorem{lemma}[theorem]{Lemma}
\newtheorem{prop}[theorem]{Proposition}
\newtheorem{corollary}[theorem]{Corollary}

{\theorembodyfont{\normalfont\rmfamily}

}

\makeatletter
\def\ack{\vspace{.5\baselineskip}\noindent{\theorem@headerfont
Acknowledgement}\ \ }
    {\renewcommand{\theththm}{\ref{#1}}\begin{ththm}}
    {\end{ththm}}
\newtheorem{ththm}{Theorem}
\newenvironment{proof}[1][]%
{\def\proof@temp{#1}\par\noindent
\textsc{Proof}\ifx\proof@temp\@empty\else\
({\proof@temp})\fi\hspace{1em}}
{\hphantom{xxx}\hfill~ {$\Box$}\par\vspace{.4\baselineskip}}
\def\operatorname#1{\mathop{\operator@font #1}\nolimits}%
\makeatother
\newcommand{\C}{\mathbb{C}}
\newcommand{\R}{\mathbb{R}}
\renewcommand{\L}{\mathcal{L}}

\newcommand{\Ad}{\operatorname{Ad}}

\newcommand{\Ker}{\operatorname{Ker}\,}

\newcommand{\half}{{\textstyle{\frac12}}}
\def\cyclic{\mathop{\kern0.9ex{{+}
\kern-2.2ex\raise-.28ex\hbox{\Large\hbox
{$\circlearrowright$}}}}\limits}

\newcommand{\egs}[1]{\arraycolsep0.1pt \renewcommand{\arraystretch}{0.5}
\begin{array}[t]{c}
=\\
#1
\end{array}
\arraycolsep5pt \renewcommand{\arraystretch}{1}}

\def\ftnote#1{\def\footnotemark{}\footnote{#1}\setcounter{footnote}{0}}

\begin{document}

\title{Construction of Ricci-type connections by reduction and induction.
\ftnote{This research was partially supported by an
Action de Recherche Concert\'ee de la Communaut\'e
fran\c{c}aise de Belgique.}\\\ }

\author{
Michel Cahen${}^{1}$
\\[20pt]
Simone Gutt${}^{1,2}$
\ftnote{\kern-4.5pt${}^{1}$Universit\'e Libre
de Bruxelles, Campus Plaine, CP 218, BE-1050~Brussels, Belgium}
\ftnote{\kern-4.5pt${}^{2}$Universit\'e de Metz, Ile du Saulcy,
F-57045~Metz Cedex 01, France}
\\[20pt]
Lorenz Schwachh\"ofer${}^{3}$
\ftnote{\kern-4.5pt${}^{3}$Mathematisches Institut, Universitaet 
Dortmund, Vogelppothsweg 87, D-44221~Dortmund, Germany}
\ftnote{Email:
\texttt{mcahen@ulb.ac.be}, \texttt{sgutt@ulb.ac.be}, 
\texttt{Lorenz.Schwachhoefer@math.uni-dortmund.de}}
 }

%\date{It is a pleasure to dedicate this paper
% to Alan Weinstein on the
%occasion of  his sixtieth birthday.}
\date{}
\setcounter{page}{0}

%% Don't want front page spaced out
\renewcommand{\baselinestretch}{1}

\maketitle

\thispagestyle{empty}

\begin{abstract}

Given the Euclidean space $\R^{2n+2}$ endowed with a constant 
symplectic structure and the standard flat connection,
and given a polynomial of degree $2$ on that space,
Baguis and Cahen \cite{bib:BaguisCahen} have defined a reduction 
procedure which yields
a symplectic manifold endowed with a Ricci-type connection. 
We observe that any symplectic manifold $(M,\omega)$
of dimension $2n ~(n\geq 2)$ endowed with a symplectic connection 
of Ricci type is locally given by a local version of such a reduction.

We also consider the reverse of this reduction procedure, an
induction procedure: we construct globally on a symplectic manifold
endowed with a connection of Ricci-type  $(M,\omega,\nabla)$ a circle
or a line bundle which embeds in a flat symplectic manifold 
$(P,\mu ,\nabla^1)$ as the
zero set of a function whose third covariant derivative vanishes,
in such a way that $(M,\omega,\nabla)$ is obtained by reduction
from $(P,\mu ,\nabla^1)$.

We further develop the particular case of symmetric symplectic manifolds
with Ricci-type connections.
 
\end{abstract}

\newpage

\renewcommand{\baselinestretch}{\linespacing}

%%%%%%%%%%%%%%%%
%%%%%%%%%%%%%%%%%
\section{Introduction}\label{sect:intro}

 A symplectic
connection $\nabla$ on a symplectic manifold $(M, \omega)$
of dimension $2n$
is a linear connection which is torsion free
and for which $\omega$ is parallel. The space of symplectic
connections on $(M,\omega)$,
${\mathcal E}(M,\omega)$ is infinite dimensional.

Selecting some particular class of connections by curvature 
conditions has, a priori, two interests. The ``moduli space''
of such particular connections may be finite dimensional; also,
on some compact symplectic manifolds which do admit a connection
of the chosen class, this connection may be ``rigid''.

In this paper, we describe completely the local behaviour
of symplectic connections of Ricci-type (see definition below)
and give some global description of simply connected 
symplectic manifolds admitting a connection of Ricci-type.

 We denote by $R$ the curvature
of $\nabla$ and by ${\underline{R}}$
the symplectic curvature tensor
$$
{\underline{R}}(X,Y,Z,T) := \omega (R(X,Y)Z,T).
$$
For any point $x \in M$, we have the symmetry properties
\begin{itemize}
\item[(i)] ${\underline{R}}_x(X,Y,Z,T)=-{\underline{R}}_x(Y,X,Z,T)$

\item[(ii)] ${\underline{R}}_x(X,Y,Z,T)={\underline{R}}_x(X,Y,T,Z)$

\item[(iii)] $\displaystyle{\cyclic_{X,Y,Z} {\underline{R}}_x(X,Y,Z,T)=0}$.
\end{itemize} 

From (i) and (ii),  $\displaystyle{{\underline{R}}_x \in
\Lambda^2 T_x^* M \otimes \odot^2 T_x^* M}$, 
where $\odot^kV$ is the  symmetrized
$k$-tensor product of the vector space $V$.

 We denote by $r$  the Ricci tensor of the connection
$\nabla$ (i.e. $r_x(X,Y) =$ tr $[Z \to R_x(X,Z)Y]$, where $X,Y,Z$ are 
in $T_xM$);
this tensor $r$ is symmetric. We denote by $\rho$ the corresponding
endomorphism of the tangent bundle:
$$
\omega(X,\rho Y) \egs{\rm def} r(X,Y)
$$ 
so that   $\rho_x$ belongs to the symplectic
algebra  $sp(T_x M, \omega_x)$; in particular tr $\rho = 0$.

The space ${\mathcal R}_x$ of symplectic curvature
tensors at $x$ is
$$
{\mathcal R}_x = {\rm ker}\, a \subset \Lambda^2 T_x^* M \otimes
\odot^2 T_x^* M
$$
where $a$ is the skewsymmetrisation map $a: \Lambda^p T_x^* M \otimes
\odot^q T_x^* M \rightarrow \Lambda^{p+1} T_x^* M \otimes
\odot^{q-1} T_x^* M$ 
$$
a(u_1 \wedge \ldots \wedge u_p \otimes v_1 \ldots v_q) :=
 \sum_{i=1}^q u_1 \wedge \ldots \wedge u_p \wedge v_i \otimes
v_1 \ldots \hat{v}_i \ldots v_q.
$$
The group $Sp(T_xM, \omega_x)$ acts on ${\mathcal R}_x$. Under
this action the space ${\mathcal R}_x,$ in dimension $2n \geq 4$, 
decomposes into two irreducible
subspaces \cite {bib:Vaisman}:
$$
{\mathcal R}_x = {\mathcal E}_x \oplus {\mathcal W}_x
$$
and the decomposition of the curvature tensor
${\underline{R}}_x$ into its ${\mathcal E}_x$ component (denoted $E_x$) and its
${\mathcal W}_x$ component (denoted $W_x$) , $
{\underline{R}}_x = E_x +W_x
$, is given by
\begin{eqnarray*}
E_x(X,Y,Z,T) &=& -\frac{1}{2(n+1)} \big[ 2 \omega_x (X,Y) r_x
(Z,T) + \omega_x(X,Z) r_x (Y,T)\\[2mm]
&\,& + \omega_x (X,T) r_x (Y,Z) - \omega_x(Y,Z) r_x(X,T)
 - \omega_x (Y,T) r_x(X,Z)\big]
\end{eqnarray*}

A connection $\nabla$ is said to be {\bf{of Ricci-type}} if,
at each point $x$, $W_x=0$.(Let us mention that such connections
were called reducible   by Vaisman in \cite{bib:Vaisman}).
In dimension 2 $(n=1)$, the space ${\mathcal W}$
vanishes identically; so we shall assume in what follows that
the manifold has dimension $m=2n>2$. 

Let us first recall two interesting features of such connections.

\noindent - When a symplectic connection is of Ricci-type,
 is satisfies the equations:
$$
\cyclic_{X,Y,Z} (\nabla_Xr)(Y,Z)=0.
$$
Those are the Euler-Lagrange equations of any natural variational principle 
whose Lagrangian  is a second degree invariant polynomial in the curvature
($r^2 $ or $  R^2$). Connections which
are solutions of those equations are called preferred ;
they are completely described in dimension $2$.

\noindent - The condition to be of Ricci-type is the condition 
on a symplectic connection $\nabla$ to have an integrable
almost complex structure $J^\nabla$ on the twistor space over $M$
 which is the bundle 
 of all compatible almost complex structures on $M$
(\cite{bib:BR}).

In this paper, we show that any symplectic manifold $(M,\omega)$
of dimension $2n ~(n\geq 2)$ admitting a symplectic connection 
of Ricci type has a local model given by a reduction procedure 
(as introduced by Baguis and Cahen in \cite{bib:BaguisCahen})
from the Euclidean space $\R^{2n+2}$ endowed with a constant 
symplectic structure and the standard flat connection.

We also consider the reverse of this reduction procedure, an
induction procedure: we construct globally on a simply connected
symplectic manifold
endowed with a connection of Ricci-type  $(M,\omega,\nabla)$ a circle
or a line bundle $N$ which embeds in a flat symplectic manifold 
$(P,\mu ,\nabla^1)$ as the
zero set of a function whose third covariant derivative vanishes,
in such a way that $(M,\omega,\nabla)$ is obtained by reduction
from $(P,\mu ,\nabla^1).$

We finally describe completely the symmetric symplectic manifolds
whose canonical connection is of Ricci-type.
Those were already studied in \cite{bib:CGRsymm} in collaboration
 with John Rawnsley.

\section{Some properties of the curvature of a Ricci-type connection}

Let $(M,\omega)$ be a smooth symplectic manifold of dim $2n$ ($n
\geq 2$) and let $\nabla$ be a smooth Ricci-type symplectic connection. 
The following results follow directly from the definition 
(and Bianchi's second identity).
 
\begin{lemma}\label{properties}
\cite{bib:CGR}
The curvature endomorphism reads
\begin{equation}\label{curv}
R(X,Y)=-\frac{1}{2(n+1)}[-2\omega (X,Y)\rho-\rho Y
\otimes \underline{X} + \rho X \otimes \underline{Y} - X \otimes
\underline{\rho Y} + Y \otimes \underline{\rho X}]
\end{equation}
 where
$\underline{X}$ denotes the 1-form $i(X) \omega$ $(for X$ a vector
field on $M$) and where, as before, $\rho$ is the endomorphism associated to the
Ricci tensor [ $ r(U,V) =\omega (U, \rho V)$ ].

Furthermore: 
\begin{itemize} 

\item[(i)] there exists a vector field $u$ such that
\begin{equation}\label{u}
\nabla_X \rho = -\dfrac{1}{2n+1} [X \otimes \underline{u} + u
\otimes \underline{X}];
\end{equation}
\item[(ii)] there exists a function $f$ such that
\begin{equation}\label{f}
\nabla_X u = -\dfrac{2n+1}{2(n+1)} \rho^2 X + f X;
\end{equation}

\item[(iii)] there exists a real number $K$ such that
\begin{equation}
tr \rho^2 + \dfrac{4(n+1)}{2n+1}f = K.
\end{equation}
\end{itemize}
\end{lemma}

\section{Construction by reduction of manifolds with Ricci type connections}
\label{section:reduction}

Let $A$ be a nonzero element in the symplectic Lie algebra
$sp(\R^{2n+2},\Omega')$ where $\Omega'$ is the 
standard symplectic structure on ${\R}^{2n+2}$.
Let $\Sigma_A$ be
the closed hypersurface $\Sigma_A \subset {\R}^{2n+2}$ with
equation :
\begin{equation}
\Omega'(x,Ax) = 1;
\end{equation}
 in order for $\Sigma_A$ to be non empty we replace, if
necessary, $A$, by $-A$.

 Let $\dot{\nabla}$ be the standard flat
symplectic affine connection on $\mathbb R^{2n+2}$. 
If $X,Y$
are vector fields tangent to $\Sigma_A$ define:
\begin{equation}
(\nabla^{\Sigma_A}_X Y)(x) = (\dot{\nabla}_X Y)(x) - \Omega'(AX,Y)x;
\end{equation}
this is a torsion free linear connection on $\Sigma_A$.

 The vector field $Ax$ is an affine vector field for this connection;
it is clearly complete and we denote by $\phi_t$ the 1-parametric group of
diffeomorphisms  of $\Sigma_A$ generated by this vector field;
clearly this flow is given by the restriction to $\Sigma_A$
of the action of $\exp tA$ on ${\R}^{2n+2}$.

Since the vector field $Ax$ is nowhere 0 on $\Sigma_A$, for any $x_0
\in \Sigma_A$, there exists :

-a neighborhood $U_{x_0}(\subset
\Sigma_A)$, 

-a ball $D\subset \R^{2n}$ of radius $r_0$, centered at the origin,

-a real interval $I=(-\epsilon,\epsilon)$ 

-and a diffeomorphism 
\begin{equation}
\chi : D \times I \to U_{x_0}
\end{equation}
 such that $\chi(0,0)=x_0$ and $\chi(y,t) =\phi_t
(\chi(y,0))$. We shall denote 
$$
\pi : U_{x_0} \to D \quad  \pi = p_1
\otimes \chi^{-1}.
$$
If we view $\Sigma_A$ as a constraint manifold
in ${\mathbb R}^{2n+2}$, $D$ is a local version of the 
Marsden-Weinstein reduction of $\Sigma_A$ around the point $x_0$. 

If $x \in \Sigma_A$, $T_x\Sigma_A = \rangle Ax \langle^{\perp}$, where
$\rangle v_1,\dots,v_p \langle$ denotes the subspace spanned by 
$v_1,\dots,v_p$ and $~^\perp$ denotes the orthogonal 
relative to $\Omega'$; let
${\mathcal H}_x (\subset T_x \Sigma_A) = \rangle x, Ax
\langle^\perp$; then
$$
T_x {\R}^{2n+2} = ({\mathcal H}_x \oplus {\R} Ax)
\oplus {\R} x
$$
and ${\pi_*}_x$ defines an isomorphism between ${\mathcal H}_x$
and the tangent space $T_yD$ for $y=\pi(x)$.
A vector belonging to ${\mathcal H}_x$ will be called horizontal.

A symplectic form on $D$,
$\omega$, is defined by
\begin{equation}
\omega_y(X,Y) = \Omega'_x (\bar{X}, \bar{Y}) \qquad y =
\pi(x)
\end{equation}
where $\bar{X}$ (resp. $\bar{Y}$) denotes the horizontal lift of
$X$ (resp. $Y$). A symplectic connection $\nabla$ on $D$ is
defined by
\begin{equation}
\overline{\nabla_X Y}(x) = \nabla^{\Sigma_A}_{\bar{X}} \bar{Y}(x) +
\Omega'(\bar{X}, \bar{Y}) Ax
\end{equation}

\begin{prop}\label{ured}
\cite{bib:BaguisCahen} The manifold $(D, \omega)$ is a symplectic manifold
and $\nabla$ is a symplectic connection of Ricci-type.

Furthermore, a direct computation shows that the corresponding $\rho, u$ and 
$f$ are given by:
\begin{eqnarray}
\overline{\rho X}(x) &=& -2(n+1)\overline{A_x}{\bar{X}}\\[2mm]
\bar{u}(x)&=&-2(n+1)(2n+1)\overline{A^2_x}x\\[2mm]
(\pi^*f)(x)&=&2(n+1)(2n+1)\Omega'(A^2x,Ax)
\end{eqnarray}
where$\overline{A^k_x}$ is the map induced by $A^k$ with values in
 ${\mathcal H}_x$:
$$
\overline{A^k_x}(X)= A^k X + \Omega'(A^k X, x) Ax - \Omega'(A^k X, Ax) x 
$$
\end{prop}

\section{Local models for symplectic connections of Ricci-type}
\label{section:localmodels}

The properties of a symplectic connection of Ricci-type, as
stated in Lemma \ref{properties}, imply in particular that

-the curvature tensor is determined by $\rho$;

-its covariant derivative is determined by $u$;

-its second covariant derivative is determined by $\rho$ and
$f$, hence by $\rho$ and $K$ with $K$ a constant;
 
-the 3rd covariant derivative of the curvature is
determined by $u, \rho, K$ and similarly for all orders.

 Hence

\begin{corollary}
Let $(M,\omega)$ be a smooth symplectic manifold of dimension $2n$ $(n
\geq 2)$ and let $\nabla$ be a smooth Ricci-type connection. Let $p_0 \in
M$; then the curvature $R_{p_0}$ and its covariant derivatives
$(\nabla^k R)_{p_0}$ (for all $k$) are determined by 
$(\rho_{x_0}, u_{x_0}, K)$.
\end{corollary}

\begin{corollary}
Let $(M, \omega, \nabla)$ (resp. $(M', \omega', \nabla')$) be two
real analytic symplectic manifolds of the same dimension $2n$ $(n \geq 2)$
each of them endowed with a symplectic connection of Ricci-type.

 Assume that there exists a linear map $b : T_{x_0}
M \to T_{x'_0}M'$ such that (i) $b^* \omega'_{x'_0} =
\omega_{x_0}$ (ii) $b u_{x_0}=u'_{x'_0}$ (iii) $b \circ \rho_{x_0}
\circ b^{-1} = \rho'_{x'_0}$. Assume further that $K =
K'$. 

Then the manifolds are locally affinely symplectically
isomorphic, i. e. there exists a normal neighborhood of $x_0$ (resp.
$x'_0$) $U_{x_0}$ (resp. $U'_{x'_0}$) and a symplectic affine
diffeomorphism $\varphi : (U_{x_0}, \omega, \nabla) \to
(U'_{x'_0}, \omega', \nabla')$ such that $\varphi(x_0) = x'_0$ and
$\varphi_{*x_0} = b$.
\end{corollary}

\noindent This follows from classical results, see for instance theorem 7.2
and corollary 7.3 in Kobayashi-Nomizu volume 1 \cite{bib:KN}.

\medskip

Consider now $(M,\omega,\nabla)$ a real
analytic symplectic manifold of dimension $2n$ $(n
\geq 2)$ endowed with an analytic Ricci-type symplectic
connection; denote as before by $ u,\rho,f$ and $K $
the associated quantities (see lemma \ref{properties}).

Let $p_0$ be a point in $M$ and choose $\xi_0$ a symplectic frame
of $T_{p_0}M$, i.e. a linear symplectic isomorphism
$\xi_0: ({\R}^{2n},\Omega)\rightarrow (T_{p_0},\omega_{p_0})$,
where $\Omega$ is the standard symplectic form on $\R^{2n}$.

Denote by $\tilde{u}(\xi_0)$ the element of $\R^{2n}$ corresponding 
to $u(p_0)$, i.e.
$$
 \tilde{u}(\xi_0)=(\xi_0)^{-1}~u(p_0)
$$
and by $\tilde{\rho}(\xi_0)$ the element of $sp(\R^{2n}, \Omega)$ 
corresponding to ${\rho}(p_0)$, i.e.
$$
\tilde{\rho}(\xi)=(\xi_0)^{-1}~\rho(p_0)~\xi_0.
$$

Define an element $A$ of $sp((\R^{2n+2}, \Omega')$ as:
$$
A= \left( \begin{array}{ccc} 0 &\dfrac{
f(p_0)}{2(n+1)(2n+1)}
&\dfrac{-{\underline{\tilde{u}(\xi_0)}}}{2(n+1)(2n+1)}\\[2mm]
1 &0 &0\\[2mm]
0 &\dfrac{- \tilde{u}(\xi_0)}{2(n+1)(2n+1)}
&\dfrac{-\tilde{\rho}(\xi_0)}{2(n+1)}
\end{array}
\right)
$$
where ${\underline{\tilde{u}(\xi_0)}}=i(\tilde{u}(\xi_0))\Omega$
and where we have chosen a basis $\{e_0,e_{0'},e_1,\dots,e_{2n}\}$
of the symplectic vector space
$R^{2n+2}$ relative to which the symplectic form
has matrix
$$
\Omega' = \left( \begin{array}{ccc} 0 &1 &0 \\-1 &0 & 0\\0&0 &\Omega
\end{array} \right)
\qquad\qquad
\Omega = \left(
\begin{array}{cc} 0 &I_n\\-I_n &0 \end{array} \right).
$$

Consider the local reduction procedure described in section 
\ref{section:reduction} from the element $A$ defined above
around the point $x_0=e_0\in
\Sigma_A=\{ x\in \R^{2n+2} ~\vert~\Omega'(x,Ax) = 1\}$.

From what we saw in section \ref{section:reduction} 
this yields a symplectic manifold with a Ricci-type connection
$(M', \omega', \nabla')$.

 Denote by $\pi'$ the map
$\pi':U_{e_0}\rightarrow M'$ where $U_{e_0}$ is the neighborhood
of ${e_0}$ in $\Sigma_A \subset \R^{2n+2}$ considered in section
\ref{section:reduction} and consider $y_0=\pi'(e_0)$.
Then ${\mathcal{H}}_{e_0}=\rangle e_0, Ae_0=e_{0'}\langle^{\perp}
=\rangle e_1, \dots, e_{2n}\langle$ is isomorphic under
$\pi'_*$ to $T_{y_0}M'$.

Introduce the injection $j : \R^{2n} \rightarrow \R^{2n+2}$
with $j(x_1,\dots, x_{2n})=(0,0,x_1,\dots, x_{2n})$ so that
$j(\R^{2n})={\mathcal{H}}_{e_0}$ and  denote by 
$b: T_{p_0}M\rightarrow T_{y_0}M'$ the map given by
$$
b=\pi'_{*_{e_0}} \circ j \circ \xi_0^{-1}.
$$
This map $b$ is a linear symplectic isomorphism since
$$
\omega'_{y_0}(bX,bY)=\Omega'(j \xi_0^{-1}X,j \xi_0^{-1}Y)=
\Omega(\xi_0^{-1}X,\xi_0^{-1}Y)=\omega_{p_0}(X,Y).
$$
Furthermore
\begin{eqnarray*}
u'(y_0)&=&\pi'_{*_{e_0}}(\bar{u'}(x_0))=
\pi'_{*_{e_0}}(-2(n+1)(2n+1)(A^2e_0-\Omega'(A^2e_0,Ae_0)e_0))\\
&=&\pi'_{*_{e_0}}(j\tilde{u}(\xi_0))=\pi'_{*_{e_0}}(\xi_0)^{-1}~u(p_0)
=bu(p_0)\\[2mm]
\rho'(y_0)bX&=&\pi'_{*_{e_0}}\overline{\rho'(y_0)X}(e_0)
=\pi'_{*_{e_0}}(-2(n+1)\overline{A_{e_0}}(j \xi_0^{-1}(X)))
=\pi'_{*_{e_0}}(j\tilde{\rho}(\xi_0)\xi_0^{-1}(X))\\
&~&{\rm{ so}}\, {\rm{ that}}\, \rho'(y_0)b=b \rho(p_0)\\[2mm]
(f')(y_0)&=&2(n+1)(2n+1)\Omega'(A^2e_0,Ae_0)=f(p_0).
\end{eqnarray*}

Hence we have
\begin{theorem}
Any real analytic symplectic manifold with a Ricci-type connection
 is locally symplectically affinely isomorphic
to the symplectic manifold with a Ricci-type connection
 obtained by a local reduction procedure
around $e_0= (1,0,\dots,0)$ from a constraint surface $\Sigma_A$
defined by a second order polynomial  
in the standard flat symplectic manifold $(\R^{2n+2},\Omega',\dot{\nabla})$.
\end{theorem}

\section{Construction of a contact manifold which is a global 
circle or line bundle over $M$} \label{section:contact}

Consider $(M,\omega,\nabla)$ a smooth symplectic manifold of dimension
$2n>2$ with a smooth Ricci-type
connection and let $B(M)\stackrel{\pi}{\to} M$ 
be the $Sp(\R^{2n},\Omega)$principal bundle of symplectic frames
over $M$. (An element in the fiber over a point $p\in M$ is 
a symplectic isomorphism $\xi:(\R^{2n},\Omega)
\rightarrow (T_pM,\omega_p)).$

As before, we consider $\tilde{u} : B(M) \to \mathbb R^{2n}$ the
$Sp(\R^{2n},\Omega)$ equivariant function given by
$$
\tilde{u}(\xi) =\xi^{-1}u(x) ~\rm{ where }~ x=\pi(\xi)
$$ 
and 
$\tilde{\rho} : B(M) \to sp(\R^{2n},\Omega)$ the $Sp(\R^{2n},\Omega)$
equivariant function given by
$$
\tilde{\rho}(\xi)=\xi^{-1}\rho(x)\xi
$$
and we define the $Sp(\R^{2n},\Omega)$ equivariant map $\tilde{A} : B(M)
\to sp(\R^{2n+2},\Omega')$
\begin{equation}
\tilde{A}(\xi) = \left( \begin{array}{ccc} 0 &\dfrac{(\pi^*
f)(\xi)}{2(n+1)(2n+1)}
&\dfrac{-{\underline{\tilde{u}(\xi)}}}{2(n+1)(2n+1)}\\[2mm]
1 &0 &0\\[2mm]
0 &\dfrac{- \tilde{u}(\xi)}{2(n+1)(2n+1)}
&\dfrac{-\tilde{\rho}(\xi)}{2(n+1)}
\end{array}
\right)
\end{equation}
where ${\underline{V}}=i(V)\Omega$ for $V$ in $\R^{2n}.$

We inject the symplectic group
$Sp(\R^{2n},\Omega)$ into $Sp({\mathbb{R}}^{2n+2},{\Omega'})$~
as the set of matrices
$$
\tilde{j}(A)=\left( \begin{array}{cc} I_2 &0\\0 &A
\end{array} \right)\qquad A\in Sp(\R^{2n},\Omega).
$$

\begin{lemma}
Define the  1-form $\alpha$ on $B(M)$, with values 
in $sp({\R}^{2n+2},{\Omega'})$  by:
\begin{equation}\label{ref:alpha}
 \quad \alpha_{\xi}(\overline{X}^{hor})  = \left(
\begin{array}{ccc} 
0 &\dfrac{-\omega_x(u,X)}{2(n+1)(2n+1)}
&\dfrac{-\underline{\widetilde{\rho(X)}(\xi)}}{2(n+1)}\\[2mm]
0 & 0 & -\underline{\tilde{X}(\xi)}\\[2mm]
\tilde{X}(\xi) &\dfrac{-\widetilde{\rho(X)}(\xi)}{2(n+1)} &0
\end{array}
\right)
\end{equation}
where $X\in T_xM$ with $x=\pi(\xi)$ and $\overline{X}^{hor}$ 
is the horizontal lift of $X$ in $T_{\xi}B(M)$, and by:
\begin{equation}
\quad \alpha(C^*) =\tilde{j}_*(C)
\end{equation}
 for all $C \in sp(\R^{2n},\Omega)$ where $C^*$ denotes
 the fundamental vertical vector field on $B(M)$ associated to $C$
 ($C^*_\xi=\frac{d}{dt}\xi.\exp tC_{\vert_0}$).

This form has the following properties:
\begin{itemize}
\item[(i)] 
$R_h^* \alpha= \Ad (\tilde{j}(h^{-1}))~\alpha
\quad\quad\forall h\in Sp(\R^{2n},\Omega);$

\item[(ii)]$
d\tilde{{A}}= -[\alpha,\tilde{A}];$

\item[(iii)]
$d\alpha +[\alpha,\alpha]=-2\tilde{A} \pi^*\omega$
\end{itemize}
\end{lemma}

When one has  a $G$-principal bundle $P\stackrel{p}{\to} M$, 
an embedding of the group $G$ in a larger group $G'$,
$j:G\rightarrow G'$,  and a $1$-form $\alpha$ with values 
in the Lie algebra
of $G'$, such that $ \alpha(C^*) =j_*(C)$
for all $C$  in the Lie algebra of $G$
and $R_h^* \alpha= \Ad ({j}(h^{-1}))~\alpha$ for
all $h$ in $G,$
one can build the $G'-$principal bundle 
$P'=P\times_G G'\stackrel{p'}{\to} M$
and the unique connection $1-$form on $P'$, $\alpha'$ satisfying
$i^*\alpha'=\alpha$ where $i:P\rightarrow P';\xi\to [(\xi,1)]$.

In our situation we build the $Sp({\mathbb{R}}^{2n+2},{\Omega'})$-
principal bundle 
$$
B'(M)=
B(M)\times_{Sp({\mathbb{R}}^{2n},{\Omega})}Sp({\mathbb{R}}^{2n+2},{\Omega'})
$$
whose elements are equivalence classes of pairs $(\xi,g)~\xi \in B(M),
g\in Sp({\mathbb{R}}^{2n+2},{\Omega'})$ with $(\xi,g)$ equivalent
to $(\xi h,\tilde{j}(h^{-1})g) ~~ \forall h\in Sp({\mathbb{R}}^{2n},{\Omega})$.

\noindent The projection $\pi':B(M)'\to M$ maps $[(\xi,g)]$ to $\pi(\xi)$.

The connection $1-$form $\alpha'$ is characterised by the fact that
$$
\alpha'_{[\xi,1]}([\overline{X}^{hor},0])=\alpha_{\xi}(\overline{X}^{hor})
$$
and the equations above give:

\begin{lemma}
The curvature $2-$form of the connection $1-$form $\alpha'$ is equal to
$ -2\tilde{A'} {\pi'}^*\omega$ where $\tilde{A'}$ is the unique
$Sp({\mathbb{R}}^{2n+2},{\Omega'})-$equivariant extension of $\tilde{A}$
to $B'(M)$.

This curvature $2-$form is invariant by parallel transport 
($d^{\alpha'}curv(\alpha')=0$).

Thus the holonomy algebra of $\alpha'$ is of dimension $1$.
\end{lemma}

\begin{corollary}
Assume $M$ is simply connected. The holonomy bundle
of $\alpha'$ is a circle or a line bundle over $M$,
$N\stackrel{\pi'}{\to} M$. This bundle has a natural
contact structure $\nu$ given by the restriction to 
$N\subset B(M)'$ of the $1-$form $-\alpha'$ (viewed as
real valued since it is valued in a $1-$dimensional algebra).
One has $d\nu=2{\pi'}^*\omega$.
\end{corollary}

It is enlightening to point out the link between the holonomy
bundle $N$ over $M$ and the constraint surface $\Sigma_A$
when one sees $M$ as obtained (locally) by reduction.
The link is only local since $\Sigma_A$ is in general not a
principal bundle over $M$; in fact in most cases the quotient of
$\Sigma_A$ by the action of the group $\exp tA$ is at best an
orbifold.

Let $A$ be a nonzero element of $sp(\R^{2n+2},\Omega')$
and let $\Sigma_A=\{ y\in \R^{2n+2}~\vert~\Omega'(y,Ay)=1\}$;
we assume as before that it is not empty.
Assume that $(M,\omega,\nabla)$ is obtained by reduction
from $\Sigma_A$ (as before, we restrict ourselves to some
open set in $\Sigma_A$).

Let $y_0$ be a point in $\Sigma_A$, let $x_0=\pi( y_0) \in M$ and
choose a symplectic frame $\xi_0$ at $x_0$.
Let $\gamma (t)$ be a curve in $M$ such that $\gamma (0)=x_0$.
Let $\xi(t)$ be the symplectic frame at $\gamma (t)$ obtained 
by parallel transport along $\gamma$ from $\xi_0$ and let 
$y(t)$ be the horizontal curve in $\Sigma_A$ lifting
$\gamma (t)$ from $y_0$ (i.e. $\pi(y (t)=\gamma (t)$ and
$\Omega'(y(t),{\dot{y}}(t))=0$).
Define the element $C(t)$ of $Sp(\R^{2n+2},\Omega')$ as the
matrix whose columns are
$$
C(t)=\left( y(t)~Ay(t)~\overline{\xi(t)}\right)
$$
where $\overline{\xi(t)}$ consists of the $2n$ vectors
which are the horizontal lifts at the point $y(t)$
of the vectors of the frame $\xi(t)$ (the image under 
the map $\xi(t)$ of the usual basis of $\R^{2n}$). Then
$$
\frac{d}{dt} C(t)_{\vert_s}=C(s). \alpha_{\xi(s)}
({\overline{{\dot{\gamma}}(s)}}^{hor})
$$
where $\alpha$ is the $1$-form on $B(M)$ defined in 
(\ref{ref:alpha}) and where ${\overline{X}}^{hor}$
is the horizontal lift of $X$ in $B(M)$; hence
${\overline{{\dot{\gamma}}(s)}}^{hor}=
{\dot{\xi}}(s)$.

Let $B'(M)$ be the $Sp(\R^{2n+2},\Omega')$-principal
bundle over $M$ considered above and let 
$[(\xi_0,\Lambda_0)]$ (where $\Lambda_0$ is an element in
$Sp(\R^{2n+2},\Omega')$) be a point of $B'(M)$ above $x_0$.
The horizontal lift of $\gamma(t)$ to $B'(M)$ starting
from $[(\xi_0,\Lambda_0)]$ lives in the holonomy subbundle
containing this point; it reads
$$
[(\xi(t),D(t))]
$$
where $\xi(t)$ has been defined above and where $D(t)$ obeys
the differential equation
$$
\frac{d}{dt} D(t)_{\vert_s}=-\alpha_{\xi(s)}({\dot{\xi}}(s)). D(s)
$$
and has initial value $\Lambda_0$.

Define the map above $\gamma$ which sends $y(t)$ to
$[(\xi(t),D(t))]$ where 
$$
D(t)=C^{-1}(t)C(0)\Lambda_0;
$$
this map sends  elements of $\Sigma_A$ to elements in the 
holonomy bundle through $[(\xi_0,\Lambda_0)]$.

The map from the holonomy bundle through $[(\xi_0,\Lambda_0)]$
to $\R^{2n+2}$ given by:
$$
[(\xi, D)] \mapsto C_0 \Lambda_0 D^{-1} e_0
$$
where $C_0$ is a fixed element in $Sp(\R^{2n+2},\Omega')$
has value in the hypersurface $\Sigma_{A'}$
where $A'=C_0 {\tilde{A}}(\xi_0) C_0^{-1}.$

\section{Embedding of the contact manifold in a flat symplectic manifold}
\label{section:flat}

Let $(M,\omega)$ be a smooth symplectic manifold of dim $2n$ ($n
\geq 2$) and let $\nabla$ be a smooth symplectic connection of Ricci-type.
Let $(N,\alpha)$ be a smooth $(2n+1)$-dimensional contact manifold 
(i.e. $\alpha$ is a smooth $1$-form such that 
$\alpha\wedge (d\alpha)^n\ne 0$ everywhere).
Let $X$ be the corresponding Reeb vector field
(i.e. $i(X)d\alpha=0$ and $\alpha(X)=1$).
Assume there exists a smooth submersion $\pi:N\rightarrow M$
such that $d\alpha=2\pi^*\omega$. Then at each point $x\in N$,
$\Ker ({\pi_*}_x)=\R X$ and $\L_X \alpha = 0$.

Remark that such a contact manifold exists always if $M$
is simply connected as we saw in the previous section.

If $U$ is a vector field on $M$ we can define its "horizontal
 lift" $\overline{U}$ on $N$ by:
 $$
 (i)~~\pi_* \overline{U}=0\quad \quad (ii) ~~\alpha(\overline{U})=0.
 $$
Let us denote by $\nu$ the $2$-form $\nu=d\alpha= 2\pi^*\omega$ on $N$. 
Define a connection ${\nabla^N}$ on $N$ by:
\begin{eqnarray*}
 {\nabla^N}_{\overline{U}}\overline{V}&=&\overline{{\nabla}_UV}-
 \nu(\overline{U},\overline{V})X  \\[2mm]
{\nabla^N}_X\overline{U}&=&{\nabla^N}_{\overline{U}} X=-\frac{1}{2(n+1)}
\overline{\rho U}
\\[2mm]
{\nabla^N}_XX&=&-\frac{1}{2(n+1)(2n+1)}\overline{u}
\end{eqnarray*}
where $\rho$ is the Ricci endomorphism of $(M,\nabla)$
and where $u$ is the vector field on $M$ appearing in 
$\nabla\rho$, see lemma \ref{properties}.
Then ${\nabla^N}$ is a torsion free connection on $N$ and 
the Reeb vector field $X$ is an affine vector field for this connection.

The curvature of this connection has the following form:
\begin{eqnarray*}
{R}^N(\overline{U},\overline{V})\overline{W} &=&
 \frac{1}{2(n+1)} [\nu(\overline{\rho V},\overline{W})\overline{U} 
 -\nu(\overline{\rho U},\overline{W})\overline{V}]\\[2mm]
{R}^N(\overline{U},\overline{V}) X &=&
\frac{1}{2(n+1)(2n+1)}[\nu(\overline{u},\overline{V})\overline{U}
-\nu(\overline{u},\overline{U})\overline{V}] \\[2mm]
{R}^N(\overline{U},X) \overline{V}&=&
\frac{1}{2(n+1)(2n+1)}\nu(\overline{u},\overline{V})\overline{U}
+ \frac{1}{2(n+1)}\nu(\overline{U},\overline{\rho V})X  \\[2mm]
{R}^N(\overline{U},X) X &=&
\frac{1}{2(n+1)(2n+1)}[-\pi^*f~ \overline{U}+
\nu(\overline{U},\overline{u})X ]
\end{eqnarray*}
where $f$ is the function appearing 
in  lemma \ref{properties}.

Consider now the embedding of the contact manifold $N$
into the symplectic manifold $(P,\mu)$ of dimension $2n+2$, 
where 
$$
P=N\times \R
$$
and, if we denote by $s$ the variable along $\R$ and let
$\theta=e^{2s}~p_1^*\alpha$ ($p_1: P \rightarrow N$), we set
$$
\mu=d\theta=2e^{2s}~ds\wedge \alpha +e^{2s}~ d\alpha
$$
and let $i:N\rightarrow P~~ x\mapsto (x,0).$
Obviously $i^*\mu=\nu$.

We now define a connection $\nabla^1$ on $P$ as follows.
If $Z$ is a vector field along $N$, we denote by the
same letter the vector field on $P$ such that
$$
(i)~~Z_{i(x)}={i_*}_xZ \quad\quad  (ii)~~[Z,\partial_s]=0.
$$
The formulas for $\nabla^1$ are:
$$
\nabla^1_Z Z'={\nabla^N}_Z Z'+\gamma(Z,Z')\partial_s
$$
where 
\begin{eqnarray*}
\gamma(Z,Z')&=&\gamma(Z',Z)        \\[2mm]
\gamma(X,X)&=&\frac{1}{2(n+1)(2n+1)}\pi^*f  \\[2mm]
\gamma(X,\overline{U})&=&-\frac{1}{2(n+1)(2n+1)}
\nu(\overline{u},\overline{U})\\[2mm]
\gamma(\overline{U},\overline{V})&=&\frac{1}{2(n+1)}
\nu(\overline{U},\overline{\rho V})
\end{eqnarray*}
and
\begin{eqnarray*}
\nabla^1_Z{\partial_s} &=&\nabla^1_{\partial_s}Z=Z \\[2mm]
\nabla^1_{\partial_s}{\partial_s}&=&\partial_s.
\end{eqnarray*}

\begin{theorem}
The connection $\nabla^1$ on $(P,\mu)$ is symplectic and has zero 
curvature.
\end{theorem}

\begin{prop}
Let $\psi(s)$ be a smooth function on $P$. Then $\psi$ has
vanishing third covariant differential if and only if
\begin{equation}\label{prop2*}
\partial^2_s\psi-2\partial_s\psi=0. 
\end{equation}
In particular the function $e^{2s}$ has this property.
\end{prop}

The procedure described above is called the induction.

Let $(P,\mu,\nabla^1)$ be as above and let $\Sigma=N$ be the
constrained submanifold defined by $e^{2s}=1.$ 
Let $Y$ be the vector field transversal to $\Sigma$
such that $i(Y)\mu=\alpha$, thus $Y=\partial_s$.
 
 Let $H$ be the $1$-parametric group generated
by $X$. Then $\Sigma/H$ can be identified with $M$
and $(M,\omega)$  is 
the classical Marsden Weinstein reduction of $(P,\mu)$ for
the constraint $\Sigma$. 

The connection $\nabla$ on $M$ is obtained from the flat 
connection $\nabla^1$ on $(P,\mu)$ by reduction.
Hence

\begin{corollary}
Any smooth simply connected symplectic manifold with a Ricci-type
connection $(M,\omega,\nabla)$ can be obtained by reduction
from an hypersurface $\Sigma$ in a
 flat symplectic manifold $(P,\mu,\nabla^1)$ defined
 by the $1-$level set of a function $\psi$ on $P$
 whose third covariant derivative vanishes.
\end{corollary}

\begin{corollary}
Any smooth simply connected symplectic manifold with a Ricci-type
connection $(M,\omega,\nabla)$ is automatically analytic.
\end{corollary}

\begin{proof}
Since $(P,\mu,\nabla^1)$ is locally symmetric, $P,\mu$ and $\nabla^1$
are real analytic and the explicit construction given 
preserves analyticity.
\end{proof}

\section{Symmetric symplectic spaces with Ricci-type connections}
\label{section:symmetric}

\begin{lemma}
The reduction construction described in section \ref{section:reduction}
 yields a locally symmetric 
symplectic space (i.e. such that the curvature tensor is parallel)
if  and only if
the element $0 \neq A \in sp(\R^{2n+2},\Omega')$ 
satisfies $A^2 =\lambda I$ for a constant $\lambda \in \R$.
\end{lemma}

\begin{proof}
Indeed the connection $\nabla$ has parallel curvature tensor
if and only if $\nabla \rho=0$ hence iff $u=0$.
From the formulas above, this is true iff
$$
\overline{A^2_x}(x)= A^2 x - \Omega'(A^2 x, Ax) x=0
$$
for any $x\in \Sigma_A$.
When $u=0$, $f$ is a constant (cf Lemma \ref{properties})
and it follows from Lemma \ref{ured} that
$\Omega'(A^2 x, Ax)$ is a constant $\lambda$.
Since $\Sigma_A$ contains a basis of $\R^{2n+2}$,
this yields $A^2 =\lambda I$. 
\end{proof}

\begin{prop}
 If $0 \neq A \in sp(\R^{2n+2},\Omega')$ 
satisfies $A^2 =\lambda I$ for a constant $\lambda \in \R$,
the quotient of $\Sigma_A$ by the action of $\exp tA$
is a manifold and the natural projection
map $\Sigma_A \to M$ is a submersion which 
endows $\Sigma_A$ with a structure of
circle or line bundle over $M$. 
\end{prop}

\begin{proof}
Consider $0 \neq A \in sp(\R^{2n+2},\Omega')$ 
so that $A^2 =\lambda I$.

-Case 1: $\lambda >0$, say $\lambda=k^2$ with $k>0$.

Then there exists a basis of $\R^{2n+2}$ in which
$$
A=\left( \begin{array}{cc} kI_{n+1} &0\\
                            0 &-kI_{n+1}
  \end{array} \right)
\qquad 
\Omega'=\left( \begin{array}{cc} 0&I_n\\
                            -I_n&0
  \end{array} \right)
$$
so that $\Sigma_A=\{ (u,v)~u,v\in\R^{n+1} ~\vert~ -2ku\cdot v=1\}$.
The flow of the vector field $Ax$ is given by $\phi_t=e^{tA}$.

\noindent The map $\pi:\Sigma_A\rightarrow TS^n=\{ (u',v')~u',v'
\in\R^{n+1}~ 
\vert u'\cdot u'=1, u'\cdot v'=0\}$ defined by
$$
\pi(u,v)=\left( \frac{u}{\Vert u\Vert},
\Vert u\Vert(v+\frac{u}{2k\Vert u\Vert^2})\right)
$$
induces a diffeomorphism between $M=\Sigma_A/_{\phi_t}$ and $TS^n$.

 $M$ is a non compact simply connected manifold
and $\Sigma_A$ is a $\R-$bundle over $TS^n$.

-Case 2 : $\lambda <0$, say $\lambda=-k^2$ with $k>0$.

One splits $V^\C$ ($V=\R^{2n+2}$) into the eigenspaces relative to $A$,
$V^\C= V_{ik}\oplus V_{-ik}$ and observe that those
subspaces are Lagrangian.
Choosing a basis $\{z_1,\dots,z_{n+1}\}$ for $V_{ik}$, consider
$\omega_{kl}:=\Omega'(z_k,\overline{z}_l)$; then
$i\omega$ is a Hermitian matrix. A change of
basis ($z'_j=\sum z_i U^i_j$) yields $\omega'=^t U\omega U$
so we can find a basis for $V_{ik}$ so that
$\omega= {-2i}I_{p,n+1-p}=-2i\left( \begin{array}{cc} I_{p} &0\\
                            0 &-I_{n+1-p}
  \end{array} \right)$.
In the basis of $\R^{2n+2}$ given by
 $e_j=\half(z_j +\overline{z}_j)\quad f_j=\frac{1}{2i}(z_j -\overline{z}_j)$
 we have:
$$
A=\left( \begin{array}{cc} 0 & -kI\\
                           kI &0
  \end{array} \right)
  \quad\quad \Omega'= \left( \begin{array}{cc} 0 &I_{p,n+1-p}\\
                           -I_{p,n+1-p} & 0
  \end{array} \right)
  $$
so that $\Sigma_A=\{ (u,v)~u,v\in\R^{n+1} ~\vert~ k\sum_{i\le p} ((u^i)^2+(v^i)^2)
-k\sum_{i> p} ((u^i)^2+(v^i)^2)=1\}$. We assume $p\ge 1$ or replace $A$ by $-A$
so that $\Sigma_A \stackrel{\sim}{=} S^{2p-1}\times \R^{2n-2p+2}$ is non empty. 
The flow  $\phi_t$ is given by the action
of $\exp tA=\left( \begin{array}{cc} \cos kt I & -\sin kt I\\
                          \sin kt I &\cos kt I
  \end{array} \right).$
  
  Then $M=\Sigma_A/_{\phi_t}=(S^{2p-1}\times \R^{2n-2p+2})/U(1)$, so this reduced manifold is:
  
  \noindent  - $M=\R^{2n}$ if $p=1$;
  
  \noindent  - $M$ is a complex line bundle of rank $q:=n+1-p$ over the complex projective
  space $P_{p-1}(\C)=S^{2p-1}/U(1)$ if $1<p\le n$;
  
  \noindent  - $M=P_{n}(\C)$ if $p=n+1$.
  
  In all those cases, $M$ is simply connected and $\Sigma_A$ is a circle bundle over $M$;
  the only compact case is $M=P_{n}(\C)$.
  
  -Case 3 : $\lambda =0$, so $A^2=0$ with $A\ne 0$. Let us denote by $p$ the rank of $A$.
  One splits $V=\R^{2n+2}$ into $V= V_0\oplus V_1 \oplus V_2$ where $V_1=Im A$ (dim $V_1=p$),
  $V_0\oplus V_1=\Ker A$ (so dim $V_0=2n+2-2p$ and $V_0$ is symplectic, since 
  $V_0\oplus V_1 = V_1^\perp$) 
  and $V_2$ is a Lagrangian subspace of $V_0^\perp$ supplementary
  to $V_1$.
  Choose a basis of $V_2$ and a corresponding basis (dual for $\Omega'$) in $V_1$
  and a symplectic basis of $V_0$ so that in those basis
  $$
   A=\left( \begin{array}{ccc} 0 & 0&0\\
                            0&0 &A'\\
                            0&0&0
  \end{array} \right)
  \quad\quad \Omega'= \left( \begin{array}{ccc} \Omega_1 &0&0\\
                                                            0 &0&I_p\\
                                                            0 &-I_p& 0
  \end{array} \right)
  $$
   and $A'$ is symmetric. Changing the basis of $V_2$ and correspondingly 
   the basis of $V_1$, one can bring $A'$ to the form $A'=I_{r,p-r}$
   so that  $\Omega'(x,Ax) = \sum_{i\le r} (w^i)^2-\sum_{r<i\le p}(w^i)^2$
   if $x=(u,v,w)$.
   
   Hence $\Sigma_A= S^{r-1}\times \R^{2n+2-r}$ if $r>1$ and $\Sigma_A$
   consists of two copies of $\R^{2n+1}$ if $r=1$.
   
   The action of $\phi_t$ on $(u,v,w)$ is given by
   $\phi_t(u,v,w)=(u,v+tA'w,w)$ so the reduced manifold is
   
   -two copies of $\R^{2n}$ (if $r=1$);
   
   -or $M=S^{r-1}\times \R^{2n+1-r}$ if $r>1$.
   
   \noindent In all cases, $M$ is a non compact manifold and $\Sigma_A$
   is a line bundle over $M$.
   \end{proof}

\begin{prop} 
  If $0 \neq A \in sp(\R^{2n+2},\Omega')$ 
satisfies $A^2 =\lambda I$ for a constant $\lambda \in \R$,
 the quotient manifold is a symmetric space and
 the connection obtained by reduction 
is the canonical symmetric connection.
\end{prop}

\begin{proof}
Any linear symplectic transformation $B$ of $\R^{2n+2}$ 
which commutes with 
$A$ obviously induces a symplectic affine transformation 
$\beta(B)$ of the reduced space $M=\Sigma_A/\phi_t$.
If $\pi$ denotes the canonical projection 
$\pi:\Sigma_A\rightarrow M$, then
$$
\beta(B)\circ \pi= \pi\circ B.
$$
In particular the symmetry at the point $x=\pi (y), y\in \Sigma_A$
is induced by
$$
B_yu=-u +2\Omega'(u,Ay)y-2\Omega'(u,y)Ay.
$$
\end{proof}

We shall now describe the tranvection group of $M$ 
(i.e. the group $G$ of affine
transformations of $M$ generated by the composition of two
symmetries).
Let us denote by $G'$ the group $G'=\{B\in Sp(R^{2n+2},\Omega')~\vert~
BA=AB\}.$
The tranvection group of $M$ 
 is clearly included in $\beta(G')$; in fact it is the 
smallest subgroup of $\beta(G')$ stable under conjugation
by a symmetry and which acts transitively on $M$.

Let $x_0=\pi(y_0)$ be a point in $M$ and let $s_{x_0}=\beta(B_{y_0})
$ be the symmetry
at this point. Consider the automorphism of $G'$
given by conjugaison by $B_{y_0}$and denote by $\sigma$
the induced automorphism of the Lie algebra ${\mathfrak{g}}'$ of $G'$.
Let ${\mathfrak{p}}'=\{ C\in {\mathfrak{g}}'~\vert ~\sigma (C)=-C\}$
and ${\mathfrak{k}}'=\{ C\in {\mathfrak{g}}'~\vert ~\sigma (C)=C\}$.

The dimension of ${\mathfrak{p}}'$ is equal to $2n$. Indeed,
in a basis $\{e_0,e_{0'},e_1,\dots,e_{2n}\}$
of $R^{2n+2}$ in which $e_0=y_0$ and $e_{0'}=Ae_{0}$ and 
$
\Omega'=\left( \begin{array}{ccc}
              0&1&0\\
              -1&0&0\\
              0&0&\Omega
              \end{array}\right),$ one has
$$              
B_{y_0}=\left( \begin{array}{ccc}
              1&0&0\\
              0&1&0\\
              0&0&-I_{2n}
              \end{array}\right)\quad\quad 
A=\left( \begin{array}{ccc}
              0&\lambda&0\\
              1&0&0\\
              0&0&A'
              \end{array}\right),
$$
${\mathfrak{g}}'=\{\left( \begin{array}{ccc}
              b&\lambda c&\underline{A'Z}\\
              c&-b&-\underline{Z}\\
              Z&A'Z&B
              \end{array}\right) \quad b,c\in\R;Z\in \R^{2n}; B\in sp(\R^{2n},\Omega)
               ~\rm{ such}~\rm{that }~Ê   BA'=A'B \},$
 and  ${\mathfrak{p}}'=\{\left( \begin{array}{ccc}
              0& 0&\underline{A'Z}\\
              0&0&-\underline{Z}\\
              Z&A'Z&0
              \end{array}\right) \quad Z\in \R^{2n} \}$.
              
Hence the Lie algebra of the transvection group is equal to
$\beta_* ({\mathfrak p}' +[{\mathfrak p}',{\mathfrak p}'])$.

\noindent In all cases the kernel of $\beta$ is given by $\exp tA$,
and the transvection group is described as follows:

 -Case 1: $\lambda >0$, say $\lambda=k^2$ with $k>0$.

In the basis of $\R^{2n+2}$ in which
$
A=\left( \begin{array}{cc} kI_{n+1} &0\\
                            0 &-kI_{n+1}
  \end{array} \right)
\qquad 
\Omega'=\left( \begin{array}{cc} 0&I_n\\
                            -I_n&0
  \end{array} \right),
$
we have
$G'=\{\left( \begin{array}{cc}
              B&0\\
              0&(^tB)^{-1}
              \end{array}\right) \quad  B\in Gl(n+1,\R) \},$
 and $\beta$ of such an element is the identity iff $B=\lambda I$
with $\lambda>0$.

The transvection group $G$ is isomorphic to $Sl({n+1},\R)$ and 
$$
TS^n=Sl(n+1,\R)/Gl(n,\R).
$$

-Case 2 : $\lambda <0$, say $\lambda=-k^2$ with $k>0$.

In the basis of $\R^{2n+2}$ in which
$A=\left( \begin{array}{cc} 0 & -kI\\
                           kI &0
  \end{array} \right)
  ~~ \Omega'= \left( \begin{array}{cc} 0 &I_{p,n+1-p}\\
                           -I_{p,n+1-p} & 0
  \end{array} \right)$
we have
$G'=\{\left( \begin{array}{cc}
              B_1&B_2\\
              -B_2&B_1
              \end{array}\right) \quad  B_1+iB_2 \in U(p,n+1-p) \},$
 and $\beta$ of such an element is the identity iff $B_1+iB_2 =\exp -ikt$.

The transvection group $G$ is isomorphic to $SU(p,n+1-p)$ and 
$$
M=SU(p,n+1-p)/U(p-1,n+1-p).
$$

 -Case 3: $\lambda =0$, rank$A=k=p+q$.
 
In the basis of $\R^{2n+2}$ in which
$
A=\left( \begin{array}{ccc} 0 &0&0\\
                            0&0&I_{pq}\\
                          0 &0&0
  \end{array} \right)$ and $
 \Omega'=\left( \begin{array}{ccc} \Omega_1&0&0\\
                                                0&0&-I\\
                                                  0&I & 0
  \end{array} \right)
  $
we have 
${\mathfrak{g}}'=\{\left( \begin{array}{ccc}
             D &0&C\\
             -^tC\Omega_1& -^tB&F\\
            0 &0&B 
  \end{array}\right) \quad  D\in sp(\R^{2n+2-2k},\omega_1),
 B\in so(p,q,\R),F \in gl(k,\R), ^tF=F,
  C\in Mat(2n+2-2k,k,\R) \}.$
  
  Then ${\mathfrak{p}}'$ is given by the elements of ${\mathfrak{g}}'$
  for which $ D=0, C=CJ, F=-JFJ, B=-JBJ$ where $J=\left( \begin{array}{cc}
              1&0\\
             0&-I_{k-1}
              \end{array}\right)$, so $C=(u\,0\,\dots\,0)$ for $u\in \R^{2n+2-2k}$,
  $F=\left( \begin{array}{cc}
            0&^tv\\
            v& 0
              \end{array}\right)$ for $v\in \R^{k-1}$
  and   $B=\left( \begin{array}{cc}
            0&^tw'\\
            w& 0
              \end{array}\right)$ for $w\in \R^{k-1}$ and $w'=I_{p-1,q}w.$
              
  Hence ${\mathfrak{p}}' \oplus [{\mathfrak{p}}',{\mathfrak{p}}']$ is the set of all
  elements in
  ${\mathfrak{g}}'$ for which  $D=0$.

 The transvection group $G$ has algebra ${\mathfrak g}$
  isomorphic to $\{(B, F, C)\}/(0,\R I_{pq},0)$
  where $B$ is any element in $ so(p,q,\R)$, $F$ is any symmetric
  real $k\times k$ matrix ,
  and $C$ is any real $(2n+2-2k)\times k$ matrix and the bracket is defined by
  $$
  [(B, F, C),(B', F', C')]=([B,B'],-^tC\Omega_1 C'+^tC'\Omega_1 C 
  -^tBF'+^tB'F ,CB'-C'B),
  $$
  so when $p+q>2$, the Levi factor is $ so(p,q,\R)$ and the radical is a $2-$step
  nilpotent algebra.
  If $p=0$ and $q=1$ the transvection group is $\R^{2n}$ and the symmetric
  space is the standard symplectic vector space.
  If $p=q=1$ or if $p=0$ and $q=2$, the transvection group is 
  solvable but not nilpotent.
  The two solvable examples are interesting for building exact
  quantisation.

{
}

\end{document}